\documentclass[a4paper,reqno]{amsart}
\usepackage{amssymb}
\usepackage{latexsym}
\usepackage{amsmath}
\usepackage{euscript}
\usepackage{graphics,color}
\usepackage[all]{xy}
\usepackage[margin=3cm]{geometry}
\usepackage{MnSymbol} 

\newcommand{\be}{\begin{equation}}
\newcommand{\ee}{\end{equation}}
\newcommand{\ba}{\begin{eqnarray}}
\newcommand{\ea}{\end{eqnarray}}
\newcommand{\baa}{\begin{eqnarray*}}
\newcommand{\eaa}{\end{eqnarray*}}
\newcommand{\bb}{\begin{thebibliography}}
\newcommand{\eb}{\end{thebibliography}}

\newcommand{\bi}[1]{\bibitem{#1}}
\newcommand{\lab}[1]{\label{#1}}
\newcommand{\re}[1]{(\ref{#1})}

\newcounter{my}
\newcommand{\he}%
   {\stepcounter{equation}\setcounter{my}%
   {\value{equation}}\setcounter{equation}0%
   }%
\newcommand{\she}%
   {\setcounter{equation}{\value{my}}%
    }%

\renewcommand\t{\tilde}

\theoremstyle{definition}

\numberwithin{equation}{section}

\begin{document}

\title[DAHA algebra of rank 1 and OPUC]{Double affine Hecke algebra of rank 1 and orthogonal polynomials on the unit circle}

\author{Satoshi Tsujimoto}
\address{Department of Applied Mathematics and Physics, Graduate School of Informatics, Kyoto University, Yoshida-Honmachi, Kyoto, Japan 606--8501}
\email{tsujimoto.satoshi.5s@kyoto-u.jp}
\author{Luc Vinet}
\address{Centre de recherches math\'ematiques
Universite de Montr\'eal, P.O. Box 6128, Centre-ville Station,
Montr\'eal (Qu\'ebec), H3C 3J7}
\email{vinet@crm.umontreal.ca}
\author{Alexei Zhedanov}
\address{Department of Mathematics, School of Information, Renmin University of China, Beijing 100872, CHINA}
\email{zhedanov@yahoo.com}

\begin{abstract}
An inifinite-dimensional representation of the double affine Hecke algebra of rank 1 and type $(C_1^{\vee},C_1)$ in which all generators are tridiagonal is presented. This representation  naturally leads to two systems of polynomials that are orthogonal on the unit circle. These polynomials can be considered as circle analogs of the Askey-Wilson polynomials. The corresponding polynomials orthogonal on an interval are constructed and discussed.
\end{abstract}


\maketitle

\section{Introduction}
This paper offers a connection between tridiagonal representations of a double affine Hecke algebra  (DAHA) and classes of orthogonal polynomials on the circle (OPUC). Let us first make comments and observations about the relevant DAHA before we state more precisely our objectives.
The rank 1 DAHA of type $(C_1^{\vee},C_1)$ is an algebra with four generators $T_i, \: i=1,2,3,4$ that have the following properties \cite{Mac}, \cite{NS}, \cite{KZ}, \cite{Sahi1}, \cite{Sahi2}, \cite{T_DAHA}:

\vspace{5mm} 

(i) The operators $T_i$ satisfy the algebraic equations
\be
(T_i - t_i)(T_i-t_i^{-1})=0, \lab{main_1} \ee 
where $t_i, \: i=1,2,3,4$ are complex numbers.

\vspace{5mm}

(ii) The product of these operators is up to a constant factor $Q$, the identity operator $\mathcal{I}$:
\be
T_1 T_2 T_3 T_4 = Q \mathcal{I}. \lab{main_2} \ee

We will assume that $t_i \ne \pm 1, \; i=1,2,3,4$. It is then easily seen that all the operators $T_i$ are invertible. This property stems from the following. Introduce the notation 
\be
\delta_i = {1 \over 2}(t_i -1/t_i), \quad \sigma_i = {1 \over 2}(t_i +1/t_i)  \lab{brev_ds} \ee
with the obvious relation
\be
\sigma_i^2-\delta_i^2=1. \lab{sd_rel} \ee
One has then 
\be
T_i = \sigma_i \mathcal{I} + \delta_i R_i \lab{T_R} \ee
and
\be
T_i^{-1} = \sigma_i \mathcal{I} - \delta_i R_i \lab{Tinv_R} \ee
where $R_i$ are involution operators satisfying the conditions
\be
R_i^2 = \mathcal{I}. \lab{R_I} \ee
Relations \re{T_R} and \re{Tinv_R} imply that we can always replace the operators $T_i$ with the reflection operators $R_i$; this will prove convenient for our purposes. In other words, we can define this DAHA as the algebra generated by four reflection operators $R_i, i=1,2,3,4$ satisfying only one condition, namely
\be
\prod_{i=1}^4 (\sigma_i \mathcal{I} + \delta_i R_i) = Q \mathcal{I}. \lab{prod_cond} \ee
A first goal of this paper is to construct a representation of this DAHA using two simple assumptions: 

\vspace{3mm}

(i) the operators $T_i, \: i=1,2,3,4$ are represented by (infinite) tridiagonal matrices;

\vspace{3mm}

(ii) the operators $T_{12}=T_1+T_2$ and $T_{34}=T_3+T_4$ 
are irreducible tridiagonal matrices in the representation.

 As will be shown these conditions are sufficient to enable the explicit construction of representations 
of all four operators $T_i,\: i=1,2,3,4$.
With this done, we will then consider the matrix pencil eigenvalue problems
\be
(R_1 - z R_2) \chi_1 =0, \quad  (R_3-z R_4) \chi_2=0, \lab{pencils} \ee 
where $z$ is a complex parameter and $\chi_1, \: \chi_2$ are generalized eigenvectors. The problem \re{pencils} 
will be seen to generate two families of OPUC. Algebraically this can be related with the  CMV matrices approach to orthogonal polynomials on the unit circle \cite{CMV}, \cite{Watkins}. We shall identify the resulting OPUC as the circle analogs of the Askey-Wilson polynomials considered in \cite{Zhe_circle}. This is the main result of our paper. 
We shall further indicate how these OPUC can be connected with two families of Askey-Wilson polynomials on the interval. In addition, we shall construct explicitly a new class of polynomials orthogonal on the interval which can be interpreted as $q$-analogs of the Bannai-Ito polynomials.

The paper is organized as follows.
In Section 2, we construct the representations of the rank 1 DAHA of type $(C_1^{\vee}, C_1)$ that the conditions (i) and (ii) entail.
In Section 3, we provide the basic properties of orthogonal polynomials on the unit circle
that are related to the DAHA representations of Section 2.
In Section 4, these OPUC are mapped to polynomials that are orthogonal on the interval $[-1, 1]$ and identified in Section 5 with the 
Askey-Wilson polynomials and with certain $q$-analogs of the Bannai-Ito polynomials.
Please note that these are different from the $q$-analogs discussed in \cite{GVZ_q_superalgebra} that arose as Racah coefficients of the quantum superalgebra $osp_q(1, 2)$. In the latter case the $q$-analogs can be obtained from the Askey-Wilson polynomials under the substitution $q$ into $-q$; this is not so for the analogs that will be of concern in this paper.
In Section 6, we consider the algebraic transformations of the DAHA generators that turn $R_2$ into a diagonal operator.
This brings additional perspectives on the relation between OPUC and DAHA representations.
In Section 7, we discuss how the central extension to the Askey-Wilson algebra appears within our framework. 
The simplest ``free'' case of the DAHA representation is considered in Section 8. It will seen to feature
 the Chebyshev  polynomials of the first, second, third and fourth kind (which are known to be special cases of the Askey-Wilson polynomials).
Finally, in Section 9 we briefly examine  the truncation conditions for the OPUC which give rise to finite-dimensional representations of the DAHA.
The paper ends with a summary of the results and concluding remarks in Section 10.

\section{Basic relations}
\setcounter{equation}{0}
We are seeking realizations of the DAHA relation \re{main_2} which can be rewritten as follows
\be
T_1 T_2 = Q T_4^{-1} T_3^{-1}. \lab{DAHA_red} \ee
Our explicit construction of the operators $T_i$ is based on the following assumptions:

\vspace{5mm}

(i) We assume that all the operators $T_i$ are tridiagonal, {i.e.} that there exists a basis $e_n, n=0,1,2,\dots$ in which
\be
T_i e_n = \xi_{n+1}^{(i)} e_{n+1} + \eta_n^{(i)} e_n + \zeta_n^{(i)} e_{n-1} , \quad n=0,1,2,\dots, \quad i=1,2,3,4 \lab{T_3diag} \ee
where $\xi_n^{(i)}, \eta_n^{(i)}, \zeta_n^{(i)}$ are some complex coefficients (it is assumed that $\zeta_0^{(i)}=0$).

\vspace{5mm}

(ii) The operators $T_{12}=T_1+T_2$ and $T_{34}=T_3+T_4$ are supposed to be irreducible. This means that the off-diagonal coefficients $\xi_n^{(12)} =\xi_n^{(1)} + \xi_n^{(2)}$ and $\zeta_n^{(12)} =\zeta_n^{(1)} + \zeta_n^{(2)}$ are nonzero for all $n=1,2,\dots$, where
\be
T_{12} e_n = \xi_{n+1}^{(12)} e_{n+1} + \eta_n^{(12)} e_n + \zeta_n^{(12)} e_{n-1}. \lab{T_12} \ee
The existence of the basis $\{e_n\}$ follows from the known fact that given a pair of generic reflection operators, say, $R_1$ and $R_2$ it is possible to construct a basis where these operators are simultaneously tridiagonal \cite{FW}, \cite{Watkins}. The requirement that the operators $T_3$ and $T_4$ of the other pair be also tridiagonal in the same basis is non-trivial and the condition for its realization will be determined in what follows.

We can extract some information about the shape of the tridiagonal operators $R_1, R_2$ with the help of the following property: if a reflection operator $R$ is tridiagonal, it decomposes into blocks of dimension one or two \cite{FW}. This means that the sequence $\xi_n^{(1)}, \: n=1,2,\dots$ (as well as the sequence $\zeta_n^{(1)}, \: n=1,2,\dots$) cannot contain two nonzero consecutive members. However if one additionally demands that the sum $R_1+R_2$ of two reflection operators is an irreducible tridiagonal operator, it is easily seen that there are no one-dimensional blocks in one of the operators, say $R_2$. The operators $R_1$ and $R_2$ should therefore look like this:
\be
 \quad R_1 =
 \begin{pmatrix}
* & * &  &    \\
  * & * &  &   \\
   &  &              * & *  \\
   &  &              * & *  \\
  &  &    & &          * & *  \\
   &  &    & &         * & *  \\
&   &  & & & & \ddots  \\
 \end{pmatrix}, \quad   R_2 =
 \begin{pmatrix}
 1 \\
& * & * &  &    \\
  &* & * &  &   \\
   & &  &               * & * \\
   & & &              * & *  \\
  &  & &   & &           * & *  \\
   &  & &   & &         * & *  \\
&   &  & & & & & \ddots  \\
 \end{pmatrix}.
\lab{M12_def} \ee 
In what follows we will assume that all entries of the  matrices $T_i, \: i=1,2,3,4$ are real. Each $2 \times 2$ block must then be of the form  
\be
B=\begin{pmatrix}
a & b    \\
  c & -a \\
 \end{pmatrix}, \quad a^2 + bc =1 . \lab{gen_2block} \ee
We can now distinguish two different possibilities for each of these $2 \times 2$ blocks:

\begin{itemize}
 \item[(i)] 
{\it $b$ and $c$ have the same sign.} Without loss of generality we can assume that $b>0$ and $c>0$. Using an appropriate  similarity transformation $B \to S^{-1} B S$ with a diagonal matrix $S$, it is possible in this case to reduce the block $B$ to
\be
B=\begin{pmatrix}
\cos \theta & \sin \theta    \\
  \sin \theta & - \cos \theta \\
 \end{pmatrix} \lab {eucl_B} \ee
with $\theta$ some real parameter. The matrix $B$ corresponds to an improper rotation in the Euclidean plane.
 \item[(ii)] 
{\it $b$ and $c$ have different signs.} Again, by an appropriate similarity transformation one can reduce the block $B$ to the canonical form    
\be
B=\begin{pmatrix}
\cosh \theta & \sinh \theta    \\
  -\sinh \theta & - \cosh \theta \\
 \end{pmatrix} \lab {pseudo_B} \ee
with some real parameter $\theta$. This corresponds to an improper rotation in the pseudo-Euclidean plane.
\end{itemize}

In general, for a given matrix $R$, the different $2 \times 2$ blocks can be either Euclidean or pseudo-Euclidean. We shall demand however that the blocks of the pairs of operators $R_1, R_2$ and $R_3, R_4$ all be of Euclidean type. In what follows we will see that this requirement imposes that the parameters $t_1,t_2,t_3,t_4$ all be imaginary.

Assuming this ``global'' Euclidean property,  we have the following realization of $R_1$ and $R_2$ (up to interchanging their roles)
\vspace{5mm}
\be
 R_1 =
 \begin{pmatrix}
  a_{0} & r_{0} &  &    \\
  r_{0} & -a_{0} &  &   \\
   &  &              a_{2} & r_{2}  \\
   &  &              r_{2} & -a_{2}  \\
  &  &    & &          a_{4} & r_{4}  \\
   &  &    & &         r_{4} & - a_{4}  \\
&   &  & & & & \ddots  \\
 \end{pmatrix}
\lab{R1_def} \ee
and
\vspace{5mm}
\be
R_2 =
 \begin{pmatrix}
1 \\
& a_{1} & r_{1} &  &    \\
  &r_{1} & -a_{1} &  &   \\
   & &  &               a_{3} & r_{3}  \\
   & & &              r_{3} & -a_{3}  \\
  &  & &   & &           a_{5} & r_{5}  \\
   &  & &   & &         r_{5} & - a_{5}  \\
&   &  & & & & & \ddots  \\
 \end{pmatrix}.
\lab{R2_def} \ee 
with
\be
r_n = \sqrt{1-a_n^2} \lab{r_def} \ee
and where $a_n$ is a parameter subjected to the restriction
\be
-1 < a_n < 1, \quad n=0,1,\dots \lab{res_a} \ee
In this realization we assume that the matrices $R_1$ and $R_2$ are both symmetric. This can always be  achieved by an appropriate similarity transformation
\be
R \to S^{-1} R S, \lab{S_sim} \ee 
where $S=diag(z_0,z_1, z_2,\dots)$ is a diagonal matrix. 

The matrices $R_3$ and $R_4$ should have the same shape as well but it is now impossible to have them symmetric. We shall assume instead that these matrices can be obtained from a pair of symmetric matrices by an appropriate similarity transformation \re{S_sim}:
\vspace{5mm}
\be
R_3 =
 \begin{pmatrix}
1 \\
& \alpha_{1} & \rho_{1}\zeta_1 &  &    \\
  &\rho_{1}{\zeta_1}^{-1} & -\alpha_{1} &  &   \\
   & &  &               \alpha_{3} & \rho_{3}\zeta_3  \\
   & & &              \rho_{3}{\zeta_3}^{-1} & -\alpha_{3}  \\
  &  & &   & &           \alpha_{5} & \rho_{5} \zeta_5 \\
   &  & &   & &         \rho_{5}{\zeta_5}^{-1} & - \alpha_{5}  \\
&   &  & & & & & \ddots  \\
 \end{pmatrix}
\lab{R3_def} \ee 
and
\vspace{5mm}
\be
 R_4 =
 \begin{pmatrix}
  \alpha_{0} & \rho_{0} \zeta_0 &  &    \\
  \rho_{0}{\zeta_0}^{-1} & -\alpha_{0} &  &   \\
   &  &              \alpha_{2} & \rho_{2}  \zeta_2\\
   &  &              \rho_{2}{\zeta_2}^{-1} & -\alpha_{2}  \\
  &  &    & &          \alpha_{4} & \rho_{4} \zeta_4 \\
   &  &    & &         \rho_{4} {\zeta_4}^{-1} & - \alpha_{4}  \\
&   &  & & & & \ddots  \\
 \end{pmatrix}
\lab{R4_def} \ee
with
\be
\rho_n= \sqrt{1-\alpha_n^2} \lab{rho_al} \ee
and where 
\be
\zeta_n = z_{n+1}/z_n \lab{zeta_z} \ee
are the coefficients connected with the action of the diagonal matrix $S$ with entries $z_0,z_1,z_2, \dots$ on the main diagonal. Note that in \cite{Jae_DAHA}, \cite{Jae} the finite-dimensional representations of the DAHA of similar shape were considered. In order to derive the coefficients $a_n, \alpha_n$ explicitly it is sufficient to note that the matrix
\be
K = T_1 T_2 - Q T_4^{-1} T_3^{-1} \lab{matr_K} \ee
should be the zero matrix if the DAHA condition \re{DAHA_red} holds. On the other hand, the matrix $K$ is a priori a five-diagonal matrix. We therefore need to set all five diagonals equal to zero.

For the outer diagonals of the matrix $K$ we have the conditions $K_{i,i+2}=K_{i+2,i}=0$ which read
\be
\delta_1 \delta_2 r_{2n} r_{2n+1} = Q \delta_3 \delta_4 \rho_{2n} \rho_{2n+1} \frac{z_{2n+2}}{z_{2n}}, \quad n=0,1,2,\dots \lab{out_1} \ee
and
\be
\delta_1 \delta_2 r_{2n} r_{2n-1} = Q \delta_3 \delta_4 \rho_{2n} \rho_{2n-1} \frac{z_{2n-1}}{z_{2n+1}}, \quad n=1,2,\dots \lab{out_2} \ee
From these equations we obtain
\be
\rho_{2n}=\gamma_0 r_{2n}{z_{2n}z_{2n+1}}, \quad \rho_{2n-1}= \gamma_1 \frac{r_{2n-1}}{ z_{2n} z_{2n-1}} \lab{rho_r} \ee
where $\gamma_0, \gamma_1$ are constants related by the restriction
\be
\delta_1 \delta_2 = Q \gamma_0 \gamma_1 \delta_3 \delta_4. \lab{res_gd} \ee
For the next diagonals, the conditions   $K_{i,i+1}=K_{i+1,i}=0$ lead to four equations:
\be
\delta_1  (\delta_2 a_{2n+1} + \sigma_2) = Q \delta_4 \gamma_0 z_{2n+1}^2  (\delta_3 \alpha_{2n+1} - \sigma_3), \lab{int_1} \ee
\be
\delta_2  (\delta_1 a_{2n} -\sigma_1) = Q \delta_3 \gamma_1 z_{2n+1}^{-2} (\delta_4 \alpha_{2n} +\sigma_4),  \lab{int_2} \ee  
\be
\delta_1  (\delta_2 a_{2n-1} - \sigma_2) = Q \delta_4 \gamma_0 z_{2n}^{2} (\delta_3 \alpha_{2n-1} + \sigma_3),  \lab{int_3} \ee  
\be
\delta_2  (\delta_1 a_{2n} +\sigma_1) = Q \delta_3 \gamma_1 z_{2n}^{-2} (\delta_4 \alpha_{2n} -\sigma_4).  \lab{int_4} \ee   
Finally, the condition that the main diagonal of equation \re{matr_K} be zero, $K_{ii}=0$, yields
\be
(-\delta_1 a_{2n} + \sigma_1)(\delta_2 a_{2n+1}+\sigma_2) = Q(\delta_4 \alpha_{2n} + \sigma_4)(-\delta_3 \alpha_{2n+1} + \sigma_3) \lab{dc_1} \ee 
and
\be
(\delta_1 a_{2n} + \sigma_1)(-\delta_2 a_{2n-1}+\sigma_2) = Q(-\delta_4 \alpha_{2n} + \sigma_4)(\delta_3 \alpha_{2n-1} + \sigma_3). \lab{dc_2} \ee
It is easy to see however, that the conditions \re{dc_1}-\re{dc_2} are not independent from the others and follow in fact from \re{int_1}-\re{int_4}.  Indeed, multiplying \re{int_1} and \re{int_2} we obtain \re{dc_1}. Similarly, multiplying \re{int_3} and \re{int_4} we recover \re{dc_2}.  
  
We can thus restrict ourselves with solving equations \re{rho_r} and \re{int_1}-\re{int_4}. The unknowns are the three sequences $a_n, \alpha_n$ and $z_n$, where $n=0,1,2, \dots$. Note that we should put 
\be
a_{-1}=\alpha_{-1}=-1 \lab{ini_aal} \ee
as initial conditions.

We initiate the solution of these equations by deriving  simple relations between $r_n$ and $\rho_n$. To that end, replace $n$ by $n-1$ in equation \re{int_1} and multiply the resulting relation by equation \re{int_3}; use then equation \re{rho_r} in order to eliminate the unknowns $z_n$. We thus get the relation
\be
\delta_3^2 (\delta_2^2 a_{2n-1}^2 - \sigma_2^2) \rho_{2n-1}^2 = \delta_2^2 (\delta_3^2 \alpha_{2n-1}^2 -\sigma_3^2) r_{2n-1}^2. \lab{al_r_1} \ee
Since $r_n, \rho_n$ and  $\sigma_i, \delta_i, \: i=1,2,3,4$ are positive, it follows from \re{al_r_1} that
\be
\delta_3 \rho_{2n-1} = \delta_2 r_{2n-1}. \lab{rr_1} \ee
Similarly, multiplying equations \re{int_2} and \re{int_4}, we arrive at the equality
\be
\delta_4 \rho_{2n} = \delta_1 r_{2n}. \lab{rr_2} \ee
Relations \re{rr_1}-\re{rr_2} together with \re{rho_r} and \re{res_gd} yield
\be
z_{2n+1} = Q z_{2n-1}, \quad z_{2n+2} = Q^{-1} z_{2n} \lab{rQr} \ee
which imply the following expressions for $z_n$:
\be
z_{2n+1}= \xi_1 Q^{n}, \quad z_{2n} = \xi_0 Q^{-n}, \quad n=0,1,2,\dots \lab{zz_expl} \ee
where $\xi_0, \xi_1$ are arbitrary nonzero constants.

Upon substituting expressions \re{zz_expl} into the system \re{int_1}-\re{int_4}, we obtain a set of linear equations for the unknowns $a_n$ and $\alpha_n$ which can be solved in an elementary way. 
From this solution we can find the constants $\gamma_0, \gamma_1$ and $\xi_1$:
\be
\gamma_0 = -\frac{t_2t_3t_4(1-t_1^2)}{Q \xi_0^2 t_1 (1-t_4^2)}, \quad \gamma_1 = -\frac{\xi_0^2 (1-t_2^2)}{t_2^2(1-t_3^2)} , \quad \xi_1 = -\frac{Q \xi_0}{t_2 t_3}. \lab{gamma_xi} \ee
The constant $\xi_0$ remains arbitrary which means that the diagonal matrix $S$ (having entries $z_n$) is defined up to an arbitrary common factor.

A conclusion that can be drawn from the solution obtained in this fashion is that the conditions $-1 < a_n <1$ and $-1 < \alpha_n <1$ corresponding to the Euclidean case, can only be realized  if all the parameters $t_1,t_2,t_3,t_4$ are imaginary.

In order to present the solution in a form suitable for further analysis, we introduce four parameters $\beta_1,\beta_2, \beta_3,\beta_4$ related to $t_1,t_2,t_3,t_4$. Assume that these parameters $\beta_i$ satisfy the conditions
\be
0<\beta_1<1, \quad 0<\beta_2<1, \quad -1<\beta_3<0, \quad -1<\beta_4<0. \lab{beta_pos} \ee 
We take the parameters $t_i$ to be given as follows in terms of the $\beta_i$ (all square roots are assumed positive) :
\be
t_1 = i \sqrt{-\beta_4/\beta_1}, \quad t_2 = i \sqrt{-\beta_1 \beta_4}, \quad  t_3 =i Q \sqrt{-  \beta_2 \beta_3}, \quad t_4 = i \sqrt{-\beta_3/\beta_2}. \lab{beta_t} \ee
Conditions \re{beta_pos} show that all parameters $t_i, \: i=1,2,3,4$ must be imaginary. We also put \cite{T_DAHA}
\be
Q= q^{-1/2} \lab{Q_q}. \ee
The solution of system \re{int_1}-\re{int_4} can now be found easily:
\be
a_{2n} = 1-2\,{\frac {\beta_{{1}} \left( 1-\beta_{{2}}\beta_{{4}}{q}^{n} \right)  \left( 1-
\beta_{{3}}\beta_{{4}}{q}^{n} \right) }{ \left( \beta_{{1}}-\beta_{{4}} \right) 
 \left( 1-g{q}^{2\,n} \right) }}, \quad a_{2n-1} = 1-2\,{\frac { \left( 1-g{q}^{n-1} \right)  \left( 1-\beta_{{1}}\beta_{{4}}{q}^
{n} \right) }{ \left( 1-\beta_{{1}}\beta_{{4}} \right)  \left( 1-g{q}^{2\,n-1}
 \right) }} \lab{a_sol} \ee
and 
\be
\alpha_{2n} = 1-2\,{\frac {\t \beta_{{1}} \left( 1-\t \beta_{{2}}\t \beta_{{4}}{q}^{n} \right)  \left( 1-
\t \beta_{{3}}\t \beta_{{4}}{q}^{n} \right) }{ \left( \t \beta_{{1}}-\t \beta_{{4}} \right) 
 \left( 1-g{q}^{2\,n} \right) }}, \quad \alpha_{2n-1} = 1-2\,{\frac { \left( 1-g{q}^{n-1} \right)  \left( 1-\t \beta_{{1}}\t \beta_{{4}}{q}^
{n} \right) }{ \left( 1-\t \beta_{{1}}\t \beta_{{4}} \right)  \left( 1-g{q}^{2\,n-1}
 \right) }} \lab{alpha_sol} \ee
where $g=\beta_1 \beta_2 \beta_3 \beta_4$ and where 
\be
\t \beta_1 = \beta_2 q^{-1/2}, \: \t \beta_2= \beta_1 q^{1/2}, \: \t \beta_3 = \beta_4  q^{1/2}, \: \t \beta_4 = \beta_3 q^{-1/2}. \lab{tb_b} \ee
Note that formulas \re{a_sol}-\re{alpha_sol} give $a_{-1}=\alpha_{-1}=-1$ which corresponds to the standard initial condition for the $L,M$ pair \cite{Simon}.

For the coefficients $r_n = \sqrt{1-a_n^2}$ and $\rho_n = \sqrt{1-\alpha_n^2}$, we have the expressions
\be
r_{2n}^2= -4\,{\frac {\beta_{{1}}\beta_{{4}} \left( 1-\beta_{{3}}\beta_{{4}}{q}^{n} \right) 
 \left( 1-\beta_{{2}}\beta_{{4}}{q}^{n} \right)  \left(1- \beta_{{1}}\beta_{{3}}{q}^{n}
 \right)  \left(1- \beta_{{1}}\beta_{{2}}{q}^{n} \right) }{ \left( \beta_{{1}}-\beta_
{{4}} \right) ^{2} \left( 1-\beta_{{1}}\beta_{{2}}\beta_{{3}}\beta_{{4}}  {q}^{
2n}   \right) ^{2}}}, \lab{r_e} \ee
\be
r_{2n-1}^2= -4\,{\frac {\beta_{{1}}\beta_{{4}}  \left( 1-{q}^{n} \right)    \left( 1-\beta_{{1}}\beta_{{4}}{q}^{n} \right) 
 \left( 1-\beta_{{2}}\beta_{{3}}{q}^{n-1} \right)  \left(1- \beta_{{1}}\beta_2 \beta_{{3}} \beta_4 {q}^{n-1}
 \right)   }{ \left( 1-\beta_{{1}}\beta_
{{4}} \right) ^{2} \left( 1-\beta_{{1}}\beta_{{2}}\beta_{{3}}\beta_{{4}}  {q}^{
2n-1}   \right) ^{2}}} \lab{r_o} \ee
and 
\be
\rho_{2n}^2 = \frac{\beta_2 \beta_3}{\beta_1 \beta_4} \: \left(\frac{\beta_1-\beta_4}{\beta_2-\beta_3}\right)^2 r_{2n}^2, \quad 
\rho_{2n-1}^2 = \frac{\beta_2 \beta_3}{\beta_1 \beta_4} \: \left(\frac{1-\beta_1\beta_4}{1-\beta_2\beta_3 q^{-1}}\right)^2 r_{2n-1}^2. \lab{rho_oe} \ee
Note that formulas \re{rho_oe} follow from \re{rr_1}-\re{rr_2} and \re{beta_t}.

\section{OPUC associated with Schur linear pencils in a DAHA}
\setcounter{equation}{0}
The operators $L$ and $M$ are closely related to polynomials orthogonal on the unit circle $\Phi_n(z)$ \cite{Simon}. Recall that these polynomials are defined through the recurrence relations
\be
\Phi_{n+1}(z) = z\Phi_n(z) - \bar {a}_n \Phi^{*}_n(z) \lab{rec_Phi} \ee
and the initial condition
\be
\Phi_0(z) =1. \lab{ini_Phi} \ee
Here 
$$
\Phi_n^*(z) = z^n \overline{\Phi_n(1/\bar z)}.
$$
The recurrence coefficients $a_n=-\overline{\Phi_{n+1}(0)}$ are usually referred to as reflection or Verblunsky parameters. Given the recurrence parameters $a_n$, the monic polynomial $\Phi_n(z) = z^n + O(z^{n-1})$ is determined uniquely through the recurrence relation \re{rec_Phi}. Under the condition $|a_n|<1$ the polynomials $\Phi_n(z)$ are orthogonal on the unit circle \cite{Simon}:
\be
\int_{\Gamma} \Phi_n(z) \overline{\Phi_{m}(z)} d \mu(z) = h_n \delta_{nm} \lab{ort_Phi} \ee
where $h_n$ are appropriate positive normalization constants, the contour $\Gamma$ is the unit circle and the measure $d\mu(z)$ is positive on $\Gamma$. If in addition, the coefficients $a_n$ are real and satisfy the condition
\be
-1<a_n<1, \lab{real_ineq} \ee
the polynomials $\Phi_n(z)=\sum_{k=0}^n A_{nk} z^k$ will have real expansion coefficients $A_{nk}$ and the orthogonality relation will become
\be
\int_0^{2\pi} \Phi_n(e^{i\theta}) \Phi_m(e^{-i\theta}) d \sigma(\theta) = h_n \: \delta_{nm}, \lab{real_ort} \ee 
with $d\sigma$ a positive and symmetric measure: $d \sigma(2\pi -\theta) = - d \sigma(\theta)$.
We shall restrict ourselves in the following to this situation where the parameters $a_n$ are real.

It is sometimes convenient to introduce another set of Laurent polynomials $\phi_n(z)$ that are related to the polynomials $\Phi_n(z)$ as follows \cite{Simon}:
\be
\phi_{2n} = z^{n} \Phi_{2n}(z^{-1}), \quad \phi_{2n+1}(z) = z^{-n} \Phi_{2n+1}(z). \lab{phi_Phi} \ee
Note that $\phi_{2n}(z)$ is a Laurent polynomial containing the monomials $z^{-n}, z^{-n+1}, \dots, z^n$ while $\phi_{2n+1}$ is a Laurent polynomial containing the monomials  $z^{-n}, z^{-n+1}, \dots, z^n, z^{n+1}$. The Laurent polynomials $\phi_n, \: n=0,1,2,\dots$ satisfy a three-term recurrence relation that can be presented in matrix form 
through the equation
\be 
(L -z M) \vec{v} =0, \lab{LMv} \ee
where $\vec{v}$ is the column vector with components $\phi_0, \phi_1, \dots, \phi_n, \dots$ and $L$ and $M$ are symmetric involution matrices having  respectively the same form as the operators $R_1$ and $R_2$ given before. The Laurent polynomials $\phi_n(z)$ are thus eigenvectors of the generalized eigenvalue problem \re{LMv}. The operator pencil $L-zM$ is sometimes called the Schur linear pencil \cite{Watkins}. Equivalently, one can present equation \re{LMv} as the ordinary eigenvalue problem
\be
U \phi_n (z) = z \phi_n(z) \lab{U_phi} \ee
with $U= ML$.
While a priori unitary, $U$ is in fact orthogonal since the operators $L$ and $M$ are symmetric:
\be
U U^{T} =\mathcal{I}, \lab{U_ort} \ee 
with $U^{T}$ denoting the transpose of $U$. The operator $U$ is 5-diagonal in the basis where the operators $L$ and $M$ look as \re{R1_def}-\re{R2_def}. This observation is in keeping with the CMV approach to orthogonal polynomials on the unit circle \cite{Simon}, \cite{Watkins}.

Let us recap our findings at this point. Starting from tridiagonal Ansatz for the DAHA operators $T_i, i=1,\dots,4$ we have arrived at  the form \re{R1_def}-\re{R2_def} of the operators $R_1, R_2$, where the recurrence coefficients $a_n$ and $\alpha_n$  have simple explicit expressions. 
Similar results were obtained for the operators $R_3, R_4$. 
The Schur pencil 
\be
(R_1 - z R_2) \vec{v}=0, \lab{LPv} \ee
is therefore seen to generate a set of polynomials $P_n(z)$ orthogonal on the unit circle
while the pencil
\be
(R_4 - z R_3) \vec{w}, =0 \lab{LPw} \ee
leads to an adjacent family $\{\t P_n(z)\}$ of OPUC . Using the explicit expressions \re{a_sol}-\re{alpha_sol} of the coefficients $a_n$ and $\alpha_n$, one can identify the corresponding polynomials $P_n(z)$ and $\t P_n(z)$ with the circle analogs of the Askey-Wilson polynomials presented in \cite{Zhe_circle} (formula (5.18)). These polynomials are orthogonal on two arcs of the unit circle.

\section{Orthogonal polynomials on the interval corresponding to OPUC}
\setcounter{equation}{0}
If all the recurrence coefficients $a_n$ are real and satisfy the positivity condition $-1<a_n<1$ it is possible to map the polynomials $\Phi_n(z)$ onto symmetric orthogonal polynomials $S_n(x)$ on the interval $[-1,1]$ \cite{DG}, \cite{Zhe_circle}. One can introduce two families of such symmetric polynomials $S_n^{(1)}(x)$ and $S_n^{(2)}(x)$ in the real variable $x=(z^{1/2}+z^{-1/2})/2$. The polynomials of the first family are obtained from the polynomials $\Phi_n(z)$ on the unit circle by the formula
\be
S_n^{(1)}(x) = \frac{2^{-n} z^{-n/2}(\Phi_n(z) + \Phi_n^*(z))}{1-a_{n-1}}.
\lab{DG_map_SP} \ee
They satisfy the recurrence relation
\be S_{n+1}^{(1)}(x) + v_n^{(1)}
S_{n-1}^{(1)}(x) = x S_n^{(1)}(x), \lab{rec_S} \ee where \be v_n^{(1)} =
\frac{1}{4} \:(1+a_{n-1})(1-a_{n-2}). \lab{v_a} \ee 
These polynomials are orthogonal on the interval $[-1,1]$
\be
\int_{-1}^1 S_n^{(1)}(x) S_m^{(1)}(x) w^{(1)}(x) dx = 0, \quad n \ne m \lab{ort_S1} \ee
where the weight function $w^{(1)}(x)$ is connected to the weight function $\rho(\theta)$ of the polynomials $\Phi_n(z)$ according to the formula:
\be
w^{(1)}(x) =
\frac{\rho(\theta)}{\sin(\theta/2)}, \quad x= \cos(\theta/2)
\lab{w_S} \ee 
with $\rho(\theta)$ given by
\be
d \sigma(\theta) = \rho(\theta) d \theta . \lab{d_sigma_rho} \ee
The second family of symmetric orthogonal polynomials $S_n^{(2)}(x)$ is defined by \cite{Zhe_circle}
\be
S_n^{(2)}(x) = \frac{2^{-n} z^{-n/2}(z\Phi_n(z) - \Phi_n^*(z))}{z-1}.
\lab{DG_map_2} \ee
They satisfy the recurrence relation 
\be S_{n+1}^{(2)}(x) + v_n^{(2)}
S_{n-1}^{(2)}(x) = x S_n^{(2)}(x), \lab{rec_S2} \ee where \be v_n^{(2)} =
\frac{1}{4} \:(1+a_{n-1})(1-a_{n}) \lab{v2_a} \ee 
and they are orthogonal on the interval $[-1,1]$
\be
\int_{-1}^1 S_n^{(2)}(x) S_m^{(2)}(x) w^{(2)}(x) dx = 0, \quad n \ne m \lab{ort_S2} \ee
with the weight function $w^{(2)}(x)$ being  \cite{Zhe_circle}
\be
w^{(2)}(x) = (1-x^2) w^{(1)}(x). \lab{w21} \ee
It is easy to show  that
\be
L_n=\frac{S_{n+1}^{(1)}(-1)}{S_n^{(1)}(-1)} = \frac{1}{2} \left(a_{n-1}-1 \right) \lab{rat_S_1} \ee
and that the polynomials $S_n^{(2)}(x)$ are obtained from the polynomials $S_n^{(1)}(x)$ by the (double) Christoffel transformation
\be
S_n^{(2)}(x) = \frac{S_{n+2}^{(1)}(x) - K_n S_n^{(1)}(x)}{x^2-1}, \lab{SS_CT} \ee
where
\be
K_n = \frac{S_{n+2}^{(1)}(1)}{S_n^{(1)}(1)} = \frac{1}{4} \left(1-a_n)(1-a_{n-1} \right). \lab{K_n} \ee
One can also construct a third family of polynomials $S_n^{(1)}(x)$ that are orthogonal on the interval, but no longer symmetric however, by performing a Christoffel transform of the polynomials $S_n^{(1)}(x)$ using the boundary point $x=-1$ \cite{DVZ}, \cite{DSVZ}:
\be
S_n^{(3)}(x) = \frac{S_{n+1}^{(1)}(x)-L_n S_n^{(1)}(x) }{x+1}, \lab{S3} \ee
with $L_n$ as in \re{rat_S_1}. The polynomials $S_n^{(3)}(x)$ satisfy the recurrence relation
\be
S_{n+1}^{(3)}(x) + \left(\frac{a_n - a_{n-1}}{2}\right) S_{n}^{(3)}(x) + \left(\frac{ 1-a_{n-1}^2 }{4}\right) S_{n-1}^{(3)}(x) = x S_n^{(3)}(x).  \lab{rec_S3} \ee
They are orthogonal on the interval $[-1,1]$
\be
\int_{-1}^1 S_n^{(3)}(x) S_m^{(3)}(x) w^{(3)}(x) dx = 0, \quad n \ne m \lab{ort_S3} \ee
with 
\be
w^{(3)}(x)= (x+1) w^{(1)}(x).  \lab{w3} \ee
Remarkably, the orthonormal polynomials $\hat S_n^{(3)}(x)$ are the eigenvectors of the symmetric Jacobi matrix $J$ that can be constructed as the sum of two involution operators \cite{DVZ} (see also \cite{CMMV2016}):
\be
J=L+M. \lab{J=L+M} \ee
Indeed, it is sufficient to verify that the recurrence coefficients in \re{rec_S3} coincide with those of the tridiagonal matrix $L+M$, where the operators $L,M$ are the same as in \re{LMv}. The diagonal entries of the matrix $J$ are $\frac{1}{2} \: (a_n-a_{n-1})$ while the off-diagonal entries are $\frac{1}{2} \: \sqrt{1-a_{n-1}^2}$. This gives exactly the recurrence relation \re{rec_S3} in their orthonormal form.

From the symmetric polynomials $S_n^{(1)}(x)$ and $S_n^{(2)}(x)$, one can construct four families of non-symmetric orthogonal polynomials $P_n^{(1,2)}(x)$ and $Q_n^{(1,2)}(x)$.
Indeed, it is well known that the even and odd orthogonal polynomials $S_n^{(1)}(x)$ can be expressed in terms of a pair of nonsymmetric orthogonal polynomials \cite{Chi}:
\be
S_{2n}^{(1)}(x) = P_n^{(1)}(x^2), \quad S_{2n+1}^{(1)}(x) = x Q_n^{(1)}(x^2) \lab{S_PQ1} \ee
and similarly
\be
S_{2n}^{(2)}(x) = P_n^{(2)}(x^2), \quad S_{2n+1}^{(2)}(x) = x Q_n^{(2)}(x^2). \lab{S_PQ2} \ee
The monic polynomials $P_n^{(1,2)}(y)$ satisfy the recurrence relation
\be
P_{n+1}^{(1,2)}(y) + B_n^{(1,2)} P_{n}^{(1,2)}(y) + U_n^{(1,2)} P_{n-1}^{(1,2)}(y) = y P_{n}^{(1,2)}(y), \lab{rec_P12} \ee
where
\ba
&&B_n^{(1)}= v_{2n}^{(1)} +  v_{2n+1}^{(1)} = \frac{1}{4} (1+a_{2n-1})(1-a_{2n-2}) + \frac{1}{4} (1+a_{2n})(1-a_{2n-1}) , \nonumber \\
&&U_n^{(1)} = v_{2n}^{(1)} v_{2n-1}^{(1)} = \frac{1}{16} (1+a_{2n-1})(1-a_{2n-2}^2)(1-a_{2n-3}) \lab{BU_1} \ea
and
\ba
&&B_n^{(2)}= v_{2n}^{(2)} +  v_{2n+1}^{(2)} = \frac{1}{2} (1+a_{2n-1})(1-a_{2n}) + \frac{1}{4} (1+a_{2n})(1-a_{2n+1}) , \nonumber \\
&&U_n^{(2)} = v_{2n}^{(2)} v_{2n-1}^{(2)} = \frac{1}{16} (1-a_{2n})(1-a_{2n-1}^2)(1+a_{2n-2}). \lab{BU_2} \ea
Similarly the monic polynomials $Q_n^{(1,2)}(x)$ satisfy the recurrence relations
\be
Q_{n+1}^{(1,2)}(y) + C_n^{(1,2)} Q_{n}^{(1,2)}(y) + V_n^{(1,2)} Q_{n-1}^{(1,2)}(y) = y Q_{n}^{(1,2)}(y), \lab{rec_Q12} \ee
where
\ba
&&C_n^{(1)}= v_{2n+2}^{(1)} +  v_{2n+1}^{(1)} = \frac{1}{4} (1+a_{2n})(1-a_{2n-1}) + \frac{1}{4} (1-a_{2n})(1+a_{2n+1}) , \nonumber \\
&&V_n^{(1)} = v_{2n}^{(1)} v_{2n+1}^{(1)} = \frac{1}{16} (1+a_{2n})(1-a_{2n-1}^2)(1-a_{2n-2}) \lab{CV_1} \ea
and
\ba
&&C_n^{(2)}= v_{2n+2}^{(2)} +  v_{2n+1}^{(2)} = \frac{1}{4} (1+a_{2n})(1-a_{2n+1}) + \frac{1}{4} (1+a_{2n+1})(1-a_{2n+2}) , \nonumber \\
&&V_n^{(2)} = v_{2n}^{(2)} v_{2n-1}^{(2)} = \frac{1}{16} (1-a_{2n+1})(1-a_{2n}^2)(1+a_{2n-1}). \lab{CV_2} \ea
The polynomials $P_n^{(1,2)}(y), Q_n^{(1,2)}(x)$ coincide (up to affine transformations of the argument) with the polynomials that are obtained under the Szeg\H{o} map from the circle to the interval (see \cite{Simon} for details of this map).

All the polynomials $P_n^{(1,2)}(x)$ and  $Q_n^{(1,2)}(x)$ are orthogonal on the interval $[0,1]$: on the one hand
\be
\int_0^1 P_n^{(1,2)}(x) P_m^{(1,2)}(x) W_{P}^{(1,2)}(x) dx = 0, \quad n\ne m \lab{ort_P12} \ee 
and, on the other
\be
\int_0^1 Q_n^{(1,2)}(x) Q_m^{(1,2)}(x) W_{Q}^{(1,2)}(x) dx = 0, \quad n\ne m . \lab{ort_Q12} \ee 
The weight functions $W_{P,Q}^{(1,2)}(x)$ that intervene here are related in a simple fashion to the weight function $w^{(1)}(x)$ of the symmetric polynomials $S_n^{(1)}(x)$. Indeed according to \cite{Chi} we have
\be
W_P^{(1)}(x) = x^{-1/2} w(x^{1/2}) \lab{W_w} \ee
and for the other weight functions we get
\be
W_Q^{(1)}(x) = x W_P^{(1)}(x), \quad W_P^{(2)}(x) = (1-x) W_P^{(1)}(x), \quad W_Q^{(2)}(x) = x (1-x) W_P^{(1)}(x). \lab{WWW} \ee
This means that the polynomials $P_n^{(2)}(x), Q_n^{(1)}(x)$ and $Q_n^{(2)}(x)$ are obtained from the polynomials $P_n^{(1)}(x)$ by the application of one or two Christoffel transformations based at the endpoints $0,1$ of the orthogonality interval.

\section{Connection to the Askey-Wilson polynomials and certain $q$-analogs of the Bannai-Ito polynomials}
\setcounter{equation}{0}

It was shown in \cite{Zhe_circle} that the polynomials $S_n(x)$ corresponding to the recurrence parameters $a_n$ \re{a_sol} are generalized Askey-Wilson polynomials. Indeed, by comparing the recurrence coefficients, it can be directly verified that the polynomials $P_n^{(1,2)}(x)$ and $Q_n^{(1,2)}(x)$ defined by \re{S_PQ1}-\re{S_PQ2} are given by
\be
P_n^{(1)}(x) = \sigma^{-n} \: V_n(\sigma x + \tau; \beta_1, \beta_2, \beta_3, \beta_4), \quad Q_n^{(1)}(x) = \sigma^{-n} \: V_n(\sigma x + \tau; q\beta_1, \beta_2, \beta_3, \beta_4) \lab{P_AW} \ee
and
\be
P_n^{(2)}(x) = \sigma^{-n} \: V_n(\sigma x + \tau; \beta_1, \beta_2, \beta_3, q\beta_4), \quad Q_n^{(2)}(x) = \sigma^{-n} \: V_n(\sigma x + \tau; q\beta_1, \beta_2, \beta_3, q\beta_4), \lab{Q_AW} \ee
where
\be
V_n(x; \beta_1, \beta_2, \beta_3, \beta_4) = \kappa_n \: {_4}\Phi_3 \left({ {q^{-n}, g q^{n-1}, \beta_1 z , \beta_1 z^{-1}  \atop \beta_1 \beta_2, \beta_1 \beta_3,\beta_1 \beta_4} ; q;q } \right), \quad x=(z+z^{-1})/2 \lab{V_AW} \ee
is the explicit expression of the Askey-Wilson polynomials depending on four parameters $\beta_1, \beta_2, \beta_3, \beta_4$ \cite{KLS}.
Here ${_4}\Phi_3$ stands for the standard basic hypergeometric function \cite{KLS} and
\be
g=\beta_1 \beta_2 \beta_3 \beta_4, \; \sigma = \frac{1}{2} \left(b_4 + b_4^{-1} - b_1 - b_{1}^{-1} \right), \; \tau= \frac{1}{2} \left(b_1 + b_1^{-1}  \right); \lab{par_P_AW} \ee
the normalization coefficient $\kappa_n$   (whose explicit expression will not be needed) ensures that the coefficient of the leading monomial $x^n$ is 1 .

For $|\beta_i|<1, \: i=1,2,3,4$, the Askey-Wilson polynomials satisfy the orthogonality condition on the interval $[-1,1]$:
\be
\int_{-1}^{1} V_n(x) V_m(x) w_{AW}(x) dx = 0, \quad   n \ne m \lab{orth_AW} \ee
where $w_{AW}(x)>0$ is the weight function of the Askey-Wilson polynomials \cite{KLS}. 

It follows from \re{P_AW}-\re{Q_AW} that the polynomials $P_n^{(1,2)}(x)$ and $Q_n^{(1,2)}(x)$ are correspondingly orthogonal on the interval $[x_1,x_2]$, where
\be
x_1 = \frac{1-\tau}{\sigma}, \quad x_2 = -\frac{1+\tau}{\sigma}. \lab{x_12} \ee
Note that due to the conditions $\beta_1 >0, \: \beta_4<0$, one has that $0<x_1<x_2<1$ so that the orthogonality interval $[x_1,x_2]$ lies inside the interval $[0,1]$.

From the orthogonality relation \re{orth_AW} and from the relations \re{S_PQ1}-\re{S_PQ1}, it follows that the symmetric polynomials $S_n^{(1,2)}(x)$ are orthogonal on two distinct symmetric intervals: $[-x_2^{1/2},-x_1^{1/2}]$ and $[x_1^{1/2},x_2^{1/2}]$, located inside the interval $[-1,1]$. In turn, this leads to the orthogonality of the corresponding OPUC $\Phi_n(z)$ on two distinct arcs of the unit circle (see \cite{Zhe_circle} for details).

We have thus related the polynomials $P_n^{(1,2)}(x)$ and $Q_n^{(1,2)}(x)$ with the ordinary Askey-Wilson polynomials. The polynomials $S_n^{(1,2)}(x)$ can hence be considered as generalized symmetric Askey-Wilson polynomials.

As in \cite{DVZ} and \cite{DSVZ}, we can also construct another family of polynomials orthogonal on the interval of the real axis. These polynomials coincide with the polynomials $S_n^{(\!3\!)}(x)$ defined by \re{S3}. They are related to the tridiagonal operators 
\be
J = R_1 + R_2, \quad \tilde J = R_3 + R_4. \lab{JJ} \ee
The corresponding eigenvalue problems 
\be
J \vec{p} = x \vec{p}, \quad \t J \vec{\t p} = x \vec{\t p} \lab{JJ_pp} \ee
generate orthonormal polynomials $\vec{p} = (\hat S_0^{(\!3\!)}, \hat S_1^{(\!3\!)}(x), \hat S_2^{(\!3\!)}(x), \dots, \hat S_n^{(\!3\!)}(x), \dots)$ and
 $\vec{\t p} = ( \hat {\t S}_0^{(\!3\!)},  \hat{\t S}_1^{(\!3\!)}(x)$, $\hat{\t S}_2^{(\!3\!)}(x), \dots,  \hat{\t S}_n^{(\!3\!)}(x), \dots)$ on the interval.

For example, the polynomials $\hat S_n(x)^{(\!3\!)}$ satisfy the three-term recurrence relation \cite{DVZ}
\begin{equation}
r_{n+1} \hat S_{n+1}^{(\!3\!)}(x) + (a_n-a_{n-1}) \hat S_n^{(\!3\!)}(x) + r_n \hat
S_{n-1}^{(\!3\!)}(x) = x \hat S_n^{(\!3\!)}(x). \lab{rec_ort_P}
\end{equation}
The monic polynomials $S_n^{(\!3\!)}(x)$ are obtained from $\hat S_n^{(\!3\!)}(x)$ by a simple renormalization. We can say that the polynomials $S_n^{(\!3\!)}(x)$ are $q$-analogs of the Bannai-Ito (BI) polynomials. 
(We wish to stress as in the introduction that these $q$-analogs of the BI polynomials are different from those in \cite{GVZ_q_superalgebra}.)
Indeed, the BI polynomials are eigenvectors of the operators $T_1+T_2$ in the degenerate version of the DAHA algebra of type $(C_1^{\vee}, C_1)$ \cite{Groenevelt}, \cite{GVZ}. Our polynomials $S_n^{(\!3\!)}(x)$ are obtained by the same Ansatz. Moreover it is possible to check directly that in the limit $q=1$ the polynomials $S_n^{(\!3\!)}(x)$ becomes the BI polynomials.

Another feature that the BI and $q$-BI polynomials share, is  their relation with the symmetric polynomials $S_n^{(1)}(x)$. Indeed we already showed on the one hand, that these polynomials are related by a single Christoffel transform. It is known on the other hand, that the BI polynomials can be obtained from the generalized (symmetric) Wilson polynomials $S_n(x)$ by a single Christoffel transform \cite{TVZ}. This corresponds to the Delsarte-Genin (DG) mapping of the unit circle \cite{DSVZ}. Our $q$-BI polynomials $S_n^{(\!3\!)}(x)$ are obtained by the same procedure from the generalized symmetric Askey-Wilson polynomials.
Finally, like the BI polynomials, the $q$-BI polynomials are orthogonal on two distinct intervals of the real axis. 

The main difference between the BI and the $q$-BI polynomials is that the BI polynomials are ``classical'', i.e. they satisfy a first order difference equation involving the reflection operator which can equivalently be presented as a second order difference equation on the classical Bannai-Ito grid. Instead, the $q$-BI polynomials $S_n^{(\!3\!)}(x)$ do not satisfy a classical equation of second order. We can nevertheless expect that the $q$-BI polynomials will satisfy a higher order difference equation.

\section{Diagonalization of one reflection operator. Direct sum of Askey-Wilson polynomials}
\setcounter{equation}{0}
The usual connection between DAHAs and orthogonal polynomials is somewhat different from the one presented here.
Instead of linear pencils, the authors in \cite{NS}, \cite{Sahi2}, \cite{KZ}, \cite{NT}, \cite{T_DAHA} consider the bilinear operators
\be
X = T_1 T_2 + T_2^{-1} T_1^{-1} \lab{X_T} \ee
and 
\be
Y= T_2 T_3 + T_3^{-1} T_2^{-1} \lab{Y_T} \ee
and observe that the eigenvalue problems associated to these operators give rise to non-symmetric Askey-Wilson polynomials \cite{NS} or direct sums of two $q$-Racah polynomials \cite{NT}.                         

We shall now show how our approach can be related to this picture. Let us first rewrite the expression for the operator $X$ in terms of reflections. Using formulas \re{T_R} and \re{Tinv_R} we have
\be
X= 2 \sigma_1 \sigma_2 + \delta_1 \delta_2 \left(R_1 R_2 + R_2 R_1 \right). \lab{X_RR} \ee
It follows that up to an affine transformation, the operator $X$ can be replaced by the simpler anticommutator of the reflection operators i.e. 
\be
X = R_1R_2+R_2R_1 .\lab{X=RR} \ee
Similarly we can take for $Y$:
\be
Y = R_2R_3+R_3R_2 .\lab{Y=RR} \ee
It is obvious that the reflection $R_2$ commutes with both operators $X$ and $Y$ (the anticommutator of two arbitrary reflections always commutes with each of the  reflection operators).
We now perform a similarity transformation that diagonalizes the reflection operator $R_2$. This means that we consider the transformed operators
\be
\t R_i = S^{-1} R_i S \lab{t_R} \ee
where $S$ is an operator such that the operator $\t R_2$ is diagonal. The operator $S$ can  easily be constructed. Indeed, since $R_2$ has the block-diagonal structure
\be
R_2 = \mbox{diag}\left([1]; M_1, M_3, \dots  \right), \lab{R_2_block} \ee
where
\be
M_i=\begin{pmatrix}
a_i & r_i    \\
  r_i & -a_i \\
 \end{pmatrix},  \lab{M_block} \ee
one can choose the operator $S$ to be of the same shape
\be
S = \mbox{diag}\left([1]; S_1, S_3, \dots  \right), \lab{Sblock} \ee
with
\be
S_i=\begin{pmatrix}
-\mu_i & \nu_i    \\
  \nu_i & \mu_i \\
 \end{pmatrix},  \quad i=1,3,5,\dots \lab{S_i} \ee
and
\be
\mu_i = \sqrt{\frac{1-a_i}{2}}, \quad \nu_i = \sqrt{\frac{1+a_i}{2}} . \lab{mu_nu} \ee
It is seen that the operator $S$ is a symmetric involution (i.e. $S^2=\mathcal{I}$) that diagonalizes the reflection $R_2$:
\be
\t R_2 = SR_2S = \mbox{diag}(1, -1, 1, -1, \dots). \lab{SRS_diag} \ee
Consider now 
what the operators X and Y become after this transformation. A simple calculation shows that the transformed $Y$ is diagonal as well:
\be
\t Y = S Y S = \mbox{diag}(y_0, y_1, y_2, \dots) \lab{diag_Y} \ee
with eigenvalues
\be
y_0 =2, \quad y_{2n-1}=y_{2n}     =   2 + \frac{4(1-q^n)(q^{-n}-gq^{-1})}{(1-\beta_1 \beta_4)(1-\beta_2 \beta_3 q^{-1})}. \lab{y_n} \ee
Note that the spectrum of the operator $Y$ coincides with the Askey-Wilson spectrum (with a double degeneration).
The operator $\t X$ also has a simple structure; it is 5-diagonal with only 3 nonzero diagonals:
\be
\t X =
 \begin{pmatrix}
B_0 & 0 & A_2 & 0 \\
0 & B_{1} & 0 & A_3 & 0    \\
 A_2 & 0 & B_2 & 0 & A_4 & 0   \\
 0  & A_3 & 0 & B_3 & 0 & A_5     \\
 0 & 0 & A_4 & 0 & B_4 & 0 \\
&   &  & & & & & \ddots  \\
 \end{pmatrix},
\lab{tX} \ee 
where the nonzero entries are
\ba
&&A_{2n} = 2 \nu_{2n-1} \mu_{2n-3} r_{2n-2} = \sqrt{(1+a_{2n-1})(1-a_{2n-2}^2)(1-a_{2n-3})}, \nonumber \\
&& A_{2n+1} = 2 \nu_{2n-1} \mu_{2n+1} r_{2n} = \sqrt{(1-a_{2n+1})(1-a_{2n}^2)(1+a_{2n-1})} \lab{A_expr} \ea
and 
\ba
&&B_{2n} = 2\left(\mu_{2n-1}^2 a_{2n}  - \nu_{2n-1}^2 a_{2n-2} \right) = a_{2n}(1-a_{2n-1}) - a_{2n-2}(1+a_{2n-1}), \nonumber \\
&&B_{2n+1} = 2\left(\mu_{2n+1}^2 a_{2n}  - \nu_{2n+1}^2 a_{2n+2} \right) = a_{2n}(1-a_{2n+1}) - a_{2n+2}(1+a_{2n+1}) . \lab{B_expr} \ea
In order to understand the structure of the eigenvectors of the operator $X$ it is convenient to write down the action of the matrix $X$ in the standard basis 
\be
\varphi_0 =\left(1,0,0, \dots   \right), \;   \varphi_1 = \left(0,1, 0, \dots  \right), \: \varphi_2 = \left(0,0, 1, 0, \dots  \right)\dots. \lab{stand_e} \ee
It is found that
\be
X \varphi_n = A_{n+2} \varphi_{n+2} + B_n \varphi_n + A_n \varphi_{n-2}. \lab{X_phi} \ee
The operator $R_2$ acts according to
\be
R_2 \varphi_n = (-1)^n \varphi_n. \lab{R_2_e} \ee
Operators X such that \re{X_phi} holds, commute
with $R_2$ for any choice of the coefficients $A_n, B_n$. This means that we can consider two independent  bases $\chi^{(0)}_n=\varphi_{2n}, \; \chi^{(1)}_n=\varphi_{2n+1}, \: n=0,1,2,\dots$ in which the matrix $X$ is a direct sum of two tridiagonal matrices $X^{(e)}$ and $X^{(o)}$ acting as follows
\be
X^{(e)} \chi_n^{(0)} = A_{2n+2} \chi_{n+1}^{(0)} + B_{2n} \chi_{n}^{(0)} + A_{2n} \chi_{n-1}^{(0)} \lab{Xe}   \ee
and 
\be
X^{(o)} \chi_n^{(1)} = A_{2n+3} \chi_{n+1}^{(1)} + B_{2n+1} \chi_{n}^{(1)} + A_{2n+1} \chi_{n-1}^{(1)}. \lab{Xo}   \ee
Introduce now two sets of monic orthogonal polynomials $P_n^{(e)}(x)$ and $P_n^{(o)}(x)$ which satisfy the recurrence relations
\be
 P_{n+1}^{(e)}(x) + B_{2n} P_n^{(e)}(x) + A_{2n}^2 P_{n-1}^{(e)}(x) = x P_n^{(e)}(x), \quad P_{-1}^{(e)}=0,\: P_0^{(e)}=1 \lab{Pe_rec} \ee
and 
\be
P_{n+1}^{(o)}(x) + B_{2n+1} P_n^{(o)}(x) + A_{2n+1}^2 P_n^{(o)}(x) = x P_n^{(o)}(x), \quad P_{-1}^{(o)}=0,\: P_0^{(o)}=1 . \lab{Po_rec} \ee
These polynomials are determined uniquely by the above recurrence relations and initial conditions. 

Let $\hat P_n^{(e)}(x)$ and $\hat P_n^{(o)}(x)$ be the corresponding orthonormal polynomials. The eigenvector $\pi =\{X_0, X_1, X_2, \dots, X_n, \dots\}$ of $X$, 
\be
X \pi = x \pi, \lab{X_pi} \ee
can then be presented as a direct sum involving two sets of orthogonal polynomials:
\be
X_{2n} = \kappa_0 \hat P_n^{(e)}(x), \; X_{2n+1} = \kappa_1 \hat P_n^{(o)}(x), \quad n=0,1,2,\dots \lab{Xcomp} \ee
where $\kappa_0$ and $\kappa_1$ are arbitrary complex parameters. Note that for a generic 5-diagonal matrix $X$, the eigenvector always depends on two arbitrary parameters.

The crucial observation is that the polynomials $P_n^{(e)}(x)$ coincide with the polynomials $P_n^{(1)}(x)$ introduced in \re{S_PQ1}, while the polynomials $P_n^{(o)}(x)$ happen to be the polynomials $Q_n^{(2)}(x)$ given in \re{S_PQ2} (strictly speaking, these polynomials coincide up to an inessential affine transformation of the argument $x \rightarrow \kappa_1 x + \kappa_0$; this difference can be eliminated by a corresponding redefinition of the operator $X$). We  thus see that the polynomials $P_n^{(e)}(x)$ and $P_n^{(o)}(x)$  both turn out to be generic  Askey-Wilson polynomials as can be observed from \re{P_AW} and \re{Q_AW}. Moreover, the polynomials $P_n^{(o)}(x)$ are obtained from the polynomials $P_n^{(e)}(x)$  by a double Christoffel transform.

This result can be compared with the observation made in \cite{KN} where the authors consider relations between polynomials orthogonal on the interval and OPUC. They introduce the Hermitian operator
\be
H= U+ U^{T} = ML+LM , \lab{H_UU} \ee
with the operators $L$ and $M$ the same symmetric block-diagonal reflection operators as in the basic eigenvalue relation \re{LMv}.
The operator $H$ commutes with both reflections $L$ and $M$. One of the reflections, say $M$, can be easily diagonalized. The operator $H$ then becomes 5-diagonal and decomposes into a direct sum of two tridiagonal operators. These two Jacobi operators generate two families of Szeg\H{o} polynomials on the interval corresponding to OPUC with recurrence parameters $a_n$. These two families are connected by a double Christoffel transform.  It is seen that the operator $X$ arising from the DAHA formalism, coincides with the operator $H$ \re{H_UU}.  In our approach the diagonalization of the operator $R_2$ is equivalent to the diagonalization of the operator $M$ in \cite{KN}. Hence we arrive at the same pair of Szeg\H{o} polynomials $P_n^{(e)}(x)$ and $P_n^{(o)}(x)$ as solutions of the eigenvalue problem for the pentadiagonal operator $X=H$. Moreover, as observed, the operator $X$ in the new basis is a direct sum of two tridiagonal operators corresponding to the above Seg\H{o} polynomials.

\section{Askey-Wilson algebra $AW(3)$}
\setcounter{equation}{0}
 The results obtained in the previous section  have a simple algebraic interpretation. Indeed, let $X=X^{(e)}\oplus X^{(o)}$ and $Y=Y^{(e)}\oplus Y^{(o)}$ with $X^{(e)}$ and $X^{(o)}$ the tridiagonal operators defined by \re{Xe} and \re{Xo}, and $Y^{(e)}$ and  $Y^{(o)}$ the diagonal operators defined by 
\be
Y^{(e)} = diag\left(y_0, y_2, y_4, \dots, y_{2n}, \dots \right), \quad Y^{(o)} = diag\left(y_1, y_3, y_5, \dots, y_{2n+1}, \dots \right). \lab{diag_Yoe} \ee
The pair of operators $Z_1=X^{(e)}$ and $Z_2=Y^{(e)}$  satisfy relations equivalent to those of the Askey-Wilson algebra \cite{ZheAW}
\be
Z_1^2 Z_2 + Z_2 Z_1^2 - (q+q^{-1})Z_1Z_2Z_1 = a_1 Z_1^2 + a_2 \{Z_1,Z_2\} + a_3 Z_1 + a_4 Z_2 +a_5 \mathcal{I} \lab{Z112} \ee
and 
\be
Z_2^2 Z_1 + Z_1 Z_2^2 - (q+q^{-1})Z_2Z_1Z_2 = a_2 Z_2^2 + a_1 \{Z_1,Z_2\} + a_3 Z_2 + a_6 Z_2 +a_7 \mathcal{I} \lab{Z221} \ee
with some structure constants $a_i, \: i=1,2,\dots, 7$.
This follows from results in \cite{ZheAW} where it was showed that the two operators: the tridiagonal $Z_1$ and the diagonal $Z_2$ realize the $AW(3)$ algebra  when $Z_1$ is the generic Jacobi matrix for the Askey-Wilson polynomials and $Z_2$ is a diagonal operator having the Askey-Wilson spectrum (grid points) as entries.

Similarly, the pair of operators $W_1=X^{(o)}$ and $W_2=Y^{(o)}$  obey the defining relations of the same Askey-Wilson algebra but with different structure constants:
\be
W_1^2 W_2 + W_2 W_1^2 - (q+q^{-1})W_1W_2W_1 = b_1 W_1^2 + b_2 \{W_1,W_2\} + b_3 W_1 + b_4 W_2 +b_5 \mathcal{I} \lab{W112} \ee
and 
\be
W_2^2 W_1 + W_1 W_2^2 - (q+q^{-1})W_2W_1W_2 = b_2 W_2^2 + b_1 \{W_1,W_2\} + b_3 W_2 + b_6 W_2 +b_7 \mathcal{I}.\lab{W221} \ee
One can combine these relations together and consider the commutation relations for the full operators $X=X^{(e)}\oplus X^{(o)}$ and $Y=Y^{(e)}\oplus Y^{(o)}$. We then arrive at the central extension of the $AW(3)$ algebra found by Koornwinder \cite{KZ}. Indeed, one can write down the relations \re{Z112} and \re{W112} as a single relation:
\ba
&&X^2 Y + Y X^2 - (q+q^{-1})XYX = (a_1 P_e+ b_1P_o) X^2 + (a_2 P_e + b_2 P_0) \{X,Y\} + \nonumber \\
&&(a_3 P_e + b_3 P_0) X + (a_4 P_e + b_4 P_o)Y +( a_5 P_e + b_5 P_o)\mathcal{I}  \lab{comb_ZW1} \ea
with $P_e$ and $P_o$ projection operators for the even and odd subspaces (i.e. projections on the spaces spanned by the bases $(e_0, e_2, e_4, \dots)$ and $(e_1,e_3,\dots)$ respectively). 
A similar formula is readily found for the other relation:
\ba
&&Y^2 X + X Y^2 - (q+q^{-1})YXY = (a_2 P_e+ b_2P_o) Y^2 + (a_1 P_e + b_1 P_0) \{X,Y\} + \nonumber \\
&&(a_3 P_e + b_3 P_0) Y + (a_6 P_e + b_6 P_o)X +( a_7 P_e + b_7 P_o)\mathcal{I}.  \lab{comb_ZW2} \ea
Clearly,
\be
P_e = \frac{1}{2} \left( \mathcal{I} + \t R_2 \right), \quad P_o = \frac{1}{2} \left( \mathcal{I} - \t R_2 \right), \lab{P_oe} \ee
where $\t R_2$ is the operator defined in \re{SRS_diag}. Note that the operator $\t R_2$ commutes with both $X$ and $Y$ and that we thus arrive at a central extension of the Askey-Wilson algebra similar to one constructed in \cite{KZ}. In this algebra the operator $\t R_2$ represents the central element. 
This algebraic approach  brings an interesting perspective on results recently obtained by Nomura and Terwilliger \cite{NT}  in the finite-dimensional case.

\section{The free case as the simplest example}
\setcounter{equation}{0}
Consider a special choice of the parameters $t_i$, namely:
\be
t_1=t_2=t_3=t_4=i. \lab{t_free} \ee 
On the one hand, according to \re{beta_t}, this corresponds to the parameters
\be
\beta_1=1, \: \beta_2 = q^{1/2}, \: \beta_3 = - q^{1/2}, \: \beta_4= -1 ; \lab{beta_free} \ee
on the other hand, in view of \re{a_sol}-\re{alpha_sol}, this leads to the recurrence parameters
\be
a_n = \alpha_n =0, \quad n=0,1,2,\dots. \lab{a_n=0} \ee
In \cite{Simon}, this choice of parameters has been called the  ``free case'' for 
 the OPUC. Indeed, in this case the polynomials are simply monomials
\be
\Phi_n(z) = z^n, \quad n=0,1,2,\dots \lab{Phi_monomial} \ee
and the orthogonality measure is the uniform Lebesgue measure on the unit circle:
\be
d \sigma(\theta) = \frac{1}{2\pi} d \theta. \lab{Leb_measure} \ee 
From \re{T_R} it is seen that the  DAHA generators coincide with the reflection operators up to a common factor:
\be
T_k= i R_k, \quad k=1,2,3,4. \lab{T_R_free} \ee 
The operators $R_1$ and $R_2$ take the very simple form 
\be
R_1 = diag\left(\sigma_1, \sigma_1, \dots, \sigma_1, \dots \right), \quad R_2 = diag\left([1], \sigma_1, \sigma_1, \dots, \sigma_1, \dots \right), \lab{R12_free} \ee
where $\sigma_1$ is the Pauli matrix
\be
\sigma_1=\begin{pmatrix}
0 & 1    \\
  1 & 0 \\
 \end{pmatrix}. \lab{Pauli_1} \ee
The operators $R_3$ and $R_4$ have simple forms as well:
\be
R_3 = diag\left([1], V_1, V_3, \dots, V_{2n-1}, \dots \right), \quad R_4 = diag\left(V_0, V_2, \dots, V_{2n}, \dots \right), \lab{R34_free} \ee
with
\be
V_{2n}=\begin{pmatrix}
0 & q^{1/2-n}    \\
  q^{n-1/2} & 0 \\
 \end{pmatrix}, \quad      V_{2n-1}=\begin{pmatrix}
0 & q^{n}    \\
  q^{-n} & 0 \\
 \end{pmatrix}.   \lab{} \ee
The DAHA relation \re{DAHA_red} now becomes
\be
R_1 R_2 = q^{-1/2} R_4 R_3 \lab{DAHA_free} \ee
as can be checked directly.

The symmetric orthogonal polynomials on the interval $S_n^{(1)}(x)$ defined by \re{DG_map_SP} become the Chebyshev polynomials of the first kind:
\be
S_n^{(1)}(x) =T_n(x), \lab{S_T} \ee 

Similarly, formula \re{DG_map_2} gives
\be
S_n^{(2)}(x) = U_n(x) \lab{S_U} \ee 
i.e. the monic Chebyshev polynomials of the second kind.

The monic Chebyshev polynomials of the first and second kind are defined by standard formulas \cite{KLS}
\be
T_0(x)=1, \quad T_n(x) = 2^{1-n} \cos\left(n\theta/2\right), \: n=1,2,\dots, \quad x = \cos(\theta/2) \lab{T_def} \ee 
and
\be
U_n(x) = 2^{-n} \frac{\sin ((n+1)\theta/2)}{\sin(\theta/2)}, \; n=0,1,2,\dots \lab{U_def} \ee
From formulas \re{v_a} and \re{v2_a}, we see that for both polynomials $T_n(x)$ and $U_n(x)$, the recurrence coefficients are the same  
\be
v_n^{(1)} = v_n^{(2)} = 1/4, \quad n=1,2,3,\dots. \lab{rec_Cheb} \ee
However, for $n=0$ there is a difference:
\be
v_0^{(1)} = 1/2, \quad v_0^{(2)} = 1/4. \lab{v_0_ch} \ee

The polynomials $S_3(x)$ in this case correspond to the Chebyshev polynomials of the third kind. Indeed, it is seen from  \re{w3} that the weight function of the polynomials $S_3(x)$ is
\be
w^{(\!3\!)}(x) = (1+x) w^{(1)}(x) = \frac{1+x}{\sqrt{1-x^2}} = \sqrt{\frac{1+x}{1-x}}. \lab{w3_cheb} \ee
This weight function is that of the Chebyshev polynomials of the third kind \cite{MH}. 
Now recall from \re{rec_ort_P} that the polynomials $S_n^{(\!3\!)}(x)$ correspond to the $q$-analogs of the Bannai-Ito polynomials. We thus see that in the ``free'' case, the $q$-analogs of the Bannai-Ito polynomials become the Chebyshev polynomials of the third kind.

\vspace{3mm} {\it Remark}. If one looks at the recurrence equation \re{rec_ort_P}, 
we see that in the ``free'' case $a_n=0$, the diagonal recurrence coefficients $b_n=a_n-a_{n-1}$ are all zero ($b_n=0$). 
Hence one might haste to conclude that
the corresponding orthogonal polynomials $S_n^{(\!3\!)}(x)$ should be symmetric. We should take note however, that $a_{-1}=-1$, and that this initial condition leads to 
\be
b_0 =-1, \: b_1=b_2=b_3 =\dots=0. \lab{b_Ch3} \ee 
The coefficient $b_0$ is therefore nonzero and as a result, the polynomials $S_n^{(\!3\!)}(x)$ are not symmetric.

Consider now the diagonalization of the operator $R_2$ described in Section 6. The operator \re{Sblock} which provides the desired transformation becomes
\be
S = \mbox{diag}\left([1]; K, K, K, \dots  \right), \lab{S_free} \ee
with
\be
K= \frac{1}{\sqrt{2}} \: \begin{pmatrix}
-1 & 1    \\
  1 & 1 \\
 \end{pmatrix}.   \lab{K_free} \ee
Thus, as in Section 6, we see that the operator $R_2$ is diagonal:
\be
\t R_2 = SR_2S = \mbox{diag}(1, -1, 1, -1, \dots). \lab{SRS_diag2} \ee
The operator $X=R_1R_2+R_2R_1$ becomes pentadiagonal
\be
\t X = S X S=
 \begin{pmatrix}
0 & 0 & \sqrt{2} & 0 & 0 &0  \\
0 & 0 & 0 & 1 & 0  & 0  \\
 \sqrt{2} & 0 & 0 & 0 & 1 & 0   \\
 0  & 1 & 0 & 0 & 0 & 1     \\
 0 & 0 & 1 & 0 & 0 & 0 & \ddots\\
0& 0 & 0 & 1 & 0 & 0  \\
&   &  & &\ddots &  &\ddots  \\
 \end{pmatrix},
\lab{tX_free} \ee 
while the operator $Y = R_2R_3 + R_3 R_2$ turns out to be diagonal also:
\be
\t Y = S Y S = \mbox{diag}(y_0, y_1, y_2, \dots) \lab{diag_Y_free} \ee
with eigenvalues
\be 
y_0=2, \quad y_{2n-1}=y_{2n}= q^n + q^{-n}, \quad n=1,2,3,\dots. \lab{y_n_free} \ee
It is easy to identify the even  and the odd orthogonal polynomials, $P^{(e)}(x)$ and $P^{(o)}(x)$ respectively, 
that correspond to the decomposition of the pentadiagonal matrix into a direct sum of two tridiagonal matrices. Indeed, from \re{Pe_rec} the polynomials  $P_n^{(e)}(x)$ are seen to satisfy the recurrence relation
\be
P_{n+1}^{(e)}(x) + (1+\delta_{n,1})P_{n-1}^{(e)}(x) = x P_n^{(e)}(x), \quad P_0^{(e)}(x)=1, \; P_{-1}^{(e)}(x)=0, \lab{Pe_rec_free} \ee
where $\delta_{n,k}$ is the Kronecker symbol, 
while the following relation is obtained from \re{Pe_rec} for the
polynomials $P_n^{(o)}$
\be
P_{n+1}^{(o)}(x) + P_{n-1}^{(o)}(x) = x P_n^{(o)}(x), \quad P_0^{(o)}(x)=1, \; P_{-1}^{(o)}(x)=0. \lab{Po_rec_free} \ee
Hence, these polynomials coincide with the Chebyshev polynomials of the first and second kind:
\be
P_n^{(e)}(x) = 2^n \: T_n(x/2), \quad P_n^{(o)}(x) = 2^n \: U_n(x/2), \quad  n=1,2,\dots . \lab{P_eo_free} \ee 
They are orthogonal on the interval $[-2,2]$ with the weight functions 
\be
w_e(x) = (4-x^2)^{-1/2}, \quad w_o(x) =(4-x^2)^{1/2}. \lab{w_eo_free} \ee
As expected, the function $w_o(x)$ is obtained from $w_e(x)$ 
by multiplying the latter with the second degree polynomial $4-x^2$ which means that the polynomials $P_n^{(o)}(x)$ are double Christoffel transforms of the polynomials $P_n^{(e)}(x)$.

The $AW$ algebraic relations \re{comb_ZW1}-\re{comb_ZW2} become very simple in the ``free'' case
\ba
&&X^2 Y + Y X^2 -(q+q^{-1}) XYX = -(q-q^{-1})^2 Y, \nonumber \\
&&Y^2 X + X Y^2 -(q+q^{-1}) YXY = -(q-q^{-1})^2 X \lab{AW_free} \ea
and they do not contain the central element $R_2$ at all. 

The Casimir operator
\be
Q_{AW} = (XY)^2 + (YX)^2 - \frac{1}{2}(q+q^{-1})^2 \left( XY^2X + YX^2Y \right) + \frac{1}{2}\left(q^2-q^{-2}\right)\left(q-q^{-1}\right) \left( X^2+Y^2 \right)  \lab{CS_free} \ee
becomes a multiple of the identity operator
\be
Q_{AW} = 2 \left(q^2-q^{-2}\right)\left(q-q^{-1}\right) \: \mathcal{I}. \lab{Q_free} \ee

\section{Finite-dimensional reduction}
\setcounter{equation}{0}
We have only considered so far infinite-dimensional representations of the DAHA.
Let us now briefly discuss some finite restrictions of these representations.
Note that all finite-dimensional representations of the DAHA were classified in \cite{OS}.
It is convenient to start with well known finite-dimensional reductions of the OPUC. Assume therefore that all the recurrence parameters $a_n$ satisfy the conditions
\be
-1 < a_n < 1 , \quad n=0,1,\dots, N-1 \lab{reg_a} \ee 
and that we have in addition
\be
a_N = \pm 1, \quad N=1,3 \dots. \lab{tr_a} \ee
It then follows that the matrices $L$ and $M$ become $N+1$-dimensional. For example, when $N=6$ one has two square matrices of size $7 \times 7$
\vspace{5mm}
\be
 R_1 =
 \begin{pmatrix}
  a_{0} & r_{0} &  &    \\
  r_{0} & -a_{0} &  &   \\
   &  &              a_{2} & r_{2}  \\
   &  &              r_{2} & -a_{2}  \\
  &  &    & &          a_{4} & r_{4}  \\
   &  &    & &         r_{4} & - a_{4}  \\
&   &  & & & & \pm1  \\
 \end{pmatrix}
\lab{R1_fin} \ee
and
\vspace{5mm}
\be
R_2 =
 \begin{pmatrix}
1 \\
& a_{1} & r_{1} &  &    \\
  &r_{1} & -a_{1} &  &   \\
   & &  &               a_{3} & r_{3}  \\
   & & &              r_{3} & -a_{3}  \\
  &  & &   & &           a_{5} & r_{5}  \\
   &  & &   & &         r_{5} & - a_{5}  \\
\end{pmatrix}.
\lab{R2_fin} \ee

When \eqref{R1_fin} and \eqref{R2_fin} are satisfied, the operator $U=ML$ is a $(N+1)$-dimensional unitary matrix. This $U$ has
$N+1$ eigenvalues $\lambda_s, \: s=0,1,2,\dots, N$ which are all distinct and belong to the unit circle \cite{Simon}, \cite{Watkins}:
\be
U \chi_s = \lambda_s \chi_s, \quad s=0,1,\dots, N, \quad \lab{eig_fin} \ee
where
\be
\lambda_s = e^{i \theta_s}, \quad 0 \le \theta_s < 2 \pi , \quad \theta_s \ne \theta_r \quad \mbox{if} \quad s \ne r. \lab{n_eq_theta} \ee
Moreover, these eigenvalues are symmetric with respect to the real line: for every eigenvalue $e^{i\theta_s} \ne \pm 1$, 
its complex conjugate $e^{-i\theta_s}$ is also an eigenvalue.
All the eigenvalues $\lambda_s$ are roots of the polynomial $\Phi_{N+1}(z)$ defined by \re{rec_Phi}, that is: 
\be
\Phi_{N+1}(\lambda_s) = 0. \lab{roots_Phi} \ee
This leads to the following finite orthogonality relation for the polynomials of degrees smaller than $N+1$ \cite{Simon}, \cite{Watkins}:
\be
\sum_{s=0}^N \rho_s \Phi_n(e^{i \theta_s}) \Phi_m(e^{-i \theta_s}) = h_n \: \delta_{nm}, \quad n,m=0,1,\dots, N \lab{fin_ort} \ee
where $h_n>0$ are normalization coefficients and where the weights $\rho_s>0$ satisfy the symmetric property:
\be
\rho_s = \rho_s' \lab{sym_rho} \ee
with $\rho_s'$ corresponding to the symmetric eigenvalue $\lambda_s'=1/\lambda_s$.
The finite-dimensional reduction of the matrices $L$ and $M$ thus leads to the finite orthogonality property \re{fin_ort} of the OPUC.

Consider now possible finite-dimensional reductions for 
our tridiagonal DAHA representations as given by 
 \re{a_sol}--\re{alpha_sol}. The truncation condition $a_N= \pm 1$ is equivalent to $r_{N}=0$ where $r_n^2=1-a_n^2$. From formulas \re{r_e}--\re{r_o} it follows that the truncation condition $a_{2N} = \pm 1$ is equivalent to one of the following relations
\be
\beta_1 \beta_2 = q^{-N}, \; \beta_1 \beta_3 = q^{-N}, \; \beta_2 \beta_4 = q^{-N}, \;  \beta_3 \beta_4 = q^{-N} \lab{ev_N_tr} \ee
while the truncation condition $a_{2N+1}= \pm 1$ is equivalent to either 
\be
\beta_1 \beta_4 = q^{-N-1}, \; \beta_2 \beta_3 = q^{-N} \;{\mbox{or}}\; \beta_1 \beta_2 \beta_3 \beta_4 = q^{-N}. \lab{od_N_tr} \ee
The last condition, $\beta_1 \beta_2 \beta_3 \beta_4 = q^{-N}$, should be taken out of consideration because it leads to singular coefficients $a_n$. It thus follows that all truncations conditions can be written in the form
\be
\beta_i \beta_k = q^{-M}, \quad M=1,2,3, \dots, \quad i \ne k, \; i,k=1,2,3,4. \lab{AW_tr} \ee
These conditions \re{AW_tr} coincide in fact with the truncation ones for the Askey-Wilson polynomials \cite{KLS}. This is not surprising because our polynomials on the unit circle are related with the Askey-Wilson polynomials on the real line.

For definiteness consider the case
\be
\beta_1 \beta_4 = q^{-N-1}, \quad N=1,2,3, \dots. \lab{b_14_tr} \ee
This corresponds to the truncation condition $a_{2N+1}=1$. It follows from \re{b_14_tr} that $\beta_1$ and $\beta_4$ should have the same sign. Without loss of generality, one can put $\beta_1>0$ and $\beta_4>0$. From formulas \re{rho_oe} and from the positivity conditions $\rho_n^2>0, \; n=0,1,\dots, 2N$, we see that the parameters $\beta_2$ and $\beta_3$ should also have the same sign and should in addition both be positive.  
This leads to the conclusion that the parameters $t_1, t_2, t_3, t_4$ should be real as per \eqref{beta_t}. Thus, in contrast to the infinite-dimensional case where the parameters $t_i$ are pure imaginary, in the finite dimensional case, these parameters should be real. One can check that the corresponding formulas for the operators $T_i, \: i=1, \dots,4$ are equivalent to those obtained in \cite{Jae_DAHA} for the finite-dimensional representations of  a DAHA of a special type (associated to the $q$-Racah polynomials).

One can consider similarly the other possible truncation conditions \re{ev_N_tr}-\re{od_N_tr}. In all cases the DAHA parameters must be real.
Note that although we have discarded the ``forbidden'' truncation condition
\be
\beta_1 \beta_2 \beta_3 \beta_4 = q^{-N}, \lab{forb_AW} \ee
one can expect that under some additional restrictions on the parameters $\beta_i$, it could be possible to construct an associated  finite-dimensional representation of the DAHA. 
Indeed it was shown in \cite{LVZ} that imposing the 
truncation condition \re{forb_AW} on the Askey-Wilson polynomials  leads,
when treated appropriately, to a new finite-dimensional system of orthogonal polynomials called the para $q$-Racah polynomials. 
Finding the OPUC and the corresponding finite-dimensional DAHA representations related to the para $q$-Racah polynomials should prove to be an interesting problem.

\section{Conclusions}
\setcounter{equation}{0}
In concluding, we would like to recap the main results that have been obtained.

\begin{enumerate}
 \item[(i)] We showed that there exists a basis where all the DAHA generators $T_i, \: i=1,\dots, 4$ are block-diagonal with all blocks of dimension $2 \times 2$ (apart from the initial entry in the operators $T_2$ and $T_3$;
 \item[(ii)] The eigenvalue problem associated to the corresponding Schur linear pencil $(R_1 - z R_2) \Phi=0$ leads to the circle analogs of the Askey-Wilson polynomials $\Phi_n(z)$. These polynomials are orthogonal on two distinct arcs of the unit circle;
 \item[(iii)] The standard eigenvalue problem $(R_1 + R_2) \psi = x \psi$ gives rise to a new family of polynomials orthogonal on the interval $[-1,1]$ of the real line which can be identified as $q$-analogs of the Bannai-Ito polynomials;
 \item[(iv)] There exists another basis in which the operator $R_2$ is diagonal. In this basis the operator $Y=R_2R_3+R_3R_2$ is diagonal with entries coinciding with the points of the Askey-Wilson grid. The operator $X=R_1 R_2+R_2R_1$ is pentadiagonal with two zero intermediate off-diagonals. This is equivalent to the statement that the operator $X$ is the direct sum $X=X^{(e)} \oplus X^{(o)}$  of two tridiagonal operators $X^{(e)}$ and $X^{(o)}$. The eigenvectors of the operator $X^{(e)}$ are expressed in terms of generic Askey-Wilson polynomials. The eigenvectors of the operator $X^{(o)}$ are given in terms of adjacent Askey-Wilson polynomials obtained by performing a double Christoffel transform on the first family of Askey-Wilson polynomials;
 \item[(v)]  The algebra generated by the operators $X$ and $Y$ is a central extension of the Askey-Wilson algebra $AW(3)$. The origin of this central extension is very simple: it stems from the decomposition of the operators $X$ and $Y$ into direct sums of more elementary operators.
 \item[(vi)]  The truncation conditions lead to finite-dimensional representations of the DAHA. There is a specific 
difference concerning the structure parameters $t_i, \: i=1,2,3,4$: they have to be real while for infinite-dimensional representations must be pure imaginary.
\end{enumerate}
There are some open problems that we wish to mention and to address in the future.
\begin{enumerate}
 \item[(i)] We have associated circle analogs of the Askey-Wilson polynomials with a linear pencil constructed from the operators $R_1$ and $R_2$. The second Schur linear pencil $R_4-z R_3$ leads to another family of the Askey-Wilson polynomials on the unit circle. What is the interpretation of this second family?  
 \item[(ii)] What are the bispectral properties of the $q$-analogs of the Bannai-Ito polynomials defined by the recurrence relation \re{rec_ort_P}? We would expect these polynomials to satisfy a Dunkl-type $q$-difference equation involving involution operators.
 \item[(iii)] 
In \cite{NS}, \cite{Sahi2}, \cite{KZ} the non-symmetric Askey-Wilson polynomials were associated with infinite-dimensional representations of the DAHA. These polynomials arise as eigenfunctions of of the operator $R_1R_2+R_2R_1$ for some realization of the operators $T_i$ in terms of q-difference operators of the Dunkl type. It would be interesting to relate these non-symmetric Askey-Wilson polynomials with the circle analogs of the Askey-Wilson polynomials presented in this paper. 
 \item[(iv)] What are the possible multivariate generalizations of our construction? 
\end{enumerate}
We hope to report on some of these questions in forthcoming publications.

\section*{Acknowledgments}
\noindent 
The authors wish to thank Paul Terwilliger for his comments on the manuscript.
The research of ST is supported by  JSPS KAKENHI (Grant Numbers 16K13761) and that of LV by 
a discovery grant of the Natural Sciences and Engineering Research Council (NSERC) of Canada.

\bb{99}

\bi{Al-Salam} W.A. Al-Salam, 
{\it Characterization theorems for orthogonal polynomials}, 
in: P. Nevai (Ed.), Orthogonal Polynomials: Theory and Practice, NATO ASI Series C: Mathematical and Physical Sciences, vol. 294, Kluwer Academic Publishers, Dordrecht,1999, pp. 1--24.

\bi{AW} R.~Askey and J.~Wilson, 
{\it Some basic hypergeometric orthogonal polynomials that generalize Jacobi polynomials}, 
Mem. Amer. Math. Soc. {\bf 54}, No. 319, (1985), 1-55.

\bi{CMV} M. J. Cantero, L. Moral and L. Vel\'azquez, 
{\it Five-diagonal matrices and zeros of orthogonal polynomials on the unit circle}, 
Lin. Alg. Appl. {\bf 362} (2003), 29--56.

\bi{CMMV2016} M.J. Cantero, F. Marcell\'an, L. Moral and L. Vel\'azquez,
{\it Darboux transformations for CMV matrices},
Adv. Math. {\bf 298} (2016), 122--206.

\bi{Chi} T. Chihara, 
{\it An Introduction to Orthogonal Polynomials}, 
Gordon and Breach, NY, 1978.

\bibitem{DG}  P. Delsarte and Y. Genin, 
{\it The split Levinson algorithm}, 
IEEE Trans. Acoust. Speech Signal Process. {\bf 34} (1986), 470--478.

\bi{DVZ} M. Derevyagin, L. Vinet and A. Zhedanov, 
{\it CMV Matrices and Little and Big -1 Jacobi Polynomials}, 
Constr.Approx. {\bf 36} (2012), 513--535, arXiv:1108.3535.

\bi{DSVZ} M. Derevyagin, S. Tsujimoto, L. Vinet and A. Zhedanov, 
{\it Bannai-Ito polynomials and dressing chains}, 
Proc. Amer. Math. Soc. {\bf 142} (2014), 4191--4206, arXiv:1211.1963.

\bi{FW} C.K. Fong and P.Y. Wu, 
{\it Band-diagonal operators}, 
Lin. Alg. Appl. {\bf 248} (1996), 185--204.

\bibitem{GVZ_q_superalgebra} V.X. Genest, L. Vinet and A. Zhedanov, 
{\it The Quantum Superalgebra ${\mathfrak {osp} _ {q}(1|2)}$ and a $q$-Generalization of the Bannai–Ito Polynomials},
Commun. Math. Phys. {\bf 344} (2016), 465--481.

\bi{GVZ} V.X. Genest, L. Vinet and A. Zhedanov, 
{\it The non-symmetric Wilson polynomials are the Bannai-Ito polynomials}, 
Proc. Amer. Math. Soc. {\bf 144} (2016), 5217--5226, arXiv:1507.02995.

\bi{GLZ} Y. A. Granovskii, I. Lutzenko and A. Zhedanov, 
{\it Mutual integrability, quadratic algebras and dynamical symmetry}, 
Ann. Physics {\bf 217} (1992), 1--20.

\bi{Groenevelt} W. Groenevelt, 
{\it Fourier transforms related to a root system of rank 1}, 
Transformation Groups {\bf 12} (1) (2007), 77--116.

\bi{Ismail} M.E.H. Ismail, 
{\it Classical and Quantum orthogonal polynomials in one variable},
Encyclopedia of Mathematics and its Applications (No. 98), Cambridge, 2005.

\bi{Jae_DAHA} Jae-Ho Lee, 
{\it Q-polynomial distance-regular graphs and a double affine Hecke algebra of rank one}, 
Lin. Alg. Appl. {\bf 439} (2013), 3184--3240, arXiv:1307.5297.

\bi{Jae} Jae-Ho Lee,
{\it Nonsymmetric Askey-Wilson polynomials and Q-polynomial distance-regular graphs}, 
J. Comb. Theory, Series A {\bf 147} (2017), 75--118, arXiv:1509.04433.

\bibitem{KN} R. Killip and I. Nenciu, 
{\it Matrix Models for Circular Ensembles}, 
IMRN {\bf 2004}, Issue 50 (2004), 2665--2701.

\bibitem{KLS} R. Koekoek, P.A. Lesky and R.F. Swarttouw, 
{\it Hypergeometric orthogonal polynomials and their $q$-analogues}.
Springer, 2010.

\bibitem{LVZ} J.-M. Lemay, L.Vinet and A.Zhedanov {\it A $q$-generalization of the para-Racah polynomials}, arXiv:1708.03368.

\bibitem{KZ} T. Koornwinder, 
{\it The Relationship between Zhedanov Algebra $AW(3)$ and the Double Affine Hecke Algebra in the Rank One Case},
SIGMA {\bf 3} (2007), 063, 15 pages.

\bibitem{Mac} I.G. Macdonald, 
{\it Affine Hecke algebras and orthogonal polynomials}, 
Cambridge University Press, 2003.

\bibitem{MH} J.C. Mason and D.C. Handscomb, 
{\it Chebyshev polynomials}, 
A CRC Press Company, 2003.

\bi{NS} M. Noumi and J.V. Stokman, 
{\it Askey-Wilson polynomials: an affine Hecke algebraic approach}, 
in Laredo Lectures on Orthogonal Polynomials and Special Functions, Nova Sci. Publ., Hauppauge, NY, 2004, 111--144, math.QA/0001033.

\bi{NT} K. Nomura, P. Terwilliger, 
{\it The universal DAHA of type $({C_1}^{\vee},C_1)$ and Leonard pairs of $q$-Racah type}, 
arXiv:1701.06089.

\bi{OS} A. Oblomkov and E. Stoica, 
{\it Finite dimensional representations of double affine Hecke algebra of rank 1}, 
J.Pure Appl. Alg. {\bf 213} (2009), 766--771, arXiv:math/0409256.

\bi{Pas} P.I. Pastro, 
{\it Orthogonal polynomials and some $q$-beta integrals of Ramanujan}, 
J. Math. Anal. Appl. {\bf 112} (1985), 517--540.

\bi{Sahi1} S. Sahi, 
{\it Some properties of Koornwinder polynomials}, in $q$-Series from a Contemporary Perspective, 
Contemp. Math. {\bf 254} (2000), 395--411.

\bi{Sahi2} S. Sahi, 
{\it Raising and lowering operators for Askey-Wilson polynomials}, 
SIGMA {\bf 3} (2007), 002, 11 pages, math.QA/0701134.

\bibitem{Simon} B. Simon, 
{\it Orthogonal Polynomials On The Unit Circle}, AMS, 2005.

\bi{T_DAHA} P. Terwilliger, 
{\it Double Affine Hecke Algebras of Rank 1 and the $\mathbb{Z}_3$-Symmetric Askey-Wilson Relations}, 
SIGMA {\bf 6} (2010), 065.

\bi{T_AW} P. Terwilliger, 
{\it The Universal Askey-Wilson Algebra}, 
SIGMA {\bf 7} (2011), 069.

\bi{TVZ} S. Tsujimoto, L. Vinet and A. Zhedanov, 
{\it Dunkl shift operators and Bannai-Ito polynomials}, 
Adv. Math. {\bf 229} (2012), 2123--2158, arXiv:1106.3512.

\bi{Watkins} D.S. Watkins, 
{\it Some perspectives on the eigenvalue problem}, 
SIAM Review {\bf 35} (1993), 430--471.

\bi{ZheAW} A. Zhedanov, 
{\it Hidden symmetry of Askey-Wilson polynomials}, 
Theoretical and Mathematical Physics {\bf 89} (1991), 1146--1157.

\bi{Zhe_circle} A. Zhedanov, 
{\it On Some Classes of Polynomials Orthogonal on Arcs of the Unit Circle Connected with Symmetric Orthogonal Polynomials on an Interval}, 
J. Approx. Theory {\bf 94} (1998), 73--106. 
\end{thebibliography}

\end{document}